\documentclass[12pt,a4]{article}
\usepackage{amsfonts}
\usepackage{amssymb}
\usepackage{epsfig}
\usepackage{psfig}
\usepackage{graphics}

\title{  On nonimbeddability of topologically trivial domains and Thin Hartogs figures of $P_2(\mathbb{C})$ into Stein spaces}
\author{Fr\'ed\'eric Sarkis}
\date{ }
\newtheorem{defi}{Definition}
\newtheorem{theo}{Theorem}
\newtheorem{prop}{Proposition}

\newtheorem{coro}{Corollary}

\newenvironment{demo}{{\em Proof.}}{}

\begin{document}
\maketitle

\begin{center}
\begin{minipage}[t]{13.5cm}
{\bf abstract.} {\footnotesize A question of Poletsky was to know
if there exists a thin Hartogs figure such that any of its
neighborhoods cannot be imbedded in Stein spaces. In
\cite{chirka}, Chirka and Ivashkovitch gave such an example
arising  in an open complex manifold. In this paper, we answer to
the question of the existence of such a figure in compact surfaces
by giving an example arising in $P_2(\mathbb{C})$. By smoothing
it, we obtain a smooth (non analytic) disc with boundary
$\overline{D} \subset P_2(\mathbb{C})$ having the same property.
Consequently, this disc intersects all algebraic curves of
$P_2(\mathbb{C})$. Moreover, as $\overline D$ is topologically
trivial, it has a neighborhood diffeomorphic to the unit ball of
$\mathbb{C}^2$. This gives a negative answer to the following
question of S. Ivashkovitch: Is the property for a domain $B$ of
$P_2(\mathbb{C})$ to be diffeormorphic to the unit ball of
$\mathbb{C}^2$ a sufficient condition for the existence of
non-constant holomorphic functions on it? }
\end{minipage}
\end{center}

{\let\thefootnote\relax

\section{Introduction}
The understanding of the relationship between the theory of
holomorphic functions on two-dimensional complex manifolds and
their differential topology has been a subject of main interest,
specially in the case of domains of $P_2(\mathbb{C})$. In
\cite{nemiro2} (see also \cite{nemiro} for a review on related
questions), Nemirovski proved that if an embedded two-sphere in
$P_2(\mathbb{C})$ is not homologous to zero, then every
holomorphic function in a neighborhood of this sphere is constant.
In \cite{sarkis} \cite{sarkis2}, we proved that if $M$ is a real
hypersurface of $P_2(\mathbb{C})$ dividing it into two domains
$\Omega_1$ and $\Omega_2$, then any holomorphic function defined
in a neighborhood of at least one of this two domains is constant.
Thus topology can be an obstruction to the existence of
non-constant holomorphic functions. In this paper, we prove that
topology is not the only obstruction. Indeed, answering a question
of S. Ivashkovitch, we prove that there exists a domain of
$P_2(\mathbb{C})$ diffeormorphic to the unit ball of
$\mathbb{C}^2$ which admits no non-constant holomorphic
function.\\
A related question is the study of the imbedding of thin Hartogs
figures in Stein manifolds. Let $\Delta$ be the unit disc of
$\mathbb{C}$, $S^1$ be its boundary, $[0,1] \subset \mathbb{C}$ be
a segment in the real line and $X$ be a complex surface. We call
thin Hartogs figure the embedding of the the set ${\cal W}=\Delta
\times \{0\}\cup S^1\times [0,1] \subset \mathbb{C}^2$ by any
continuous imbedding $f:{\cal W} \rightarrow X$ which is
holomorphic on $\Delta \times \{0\}$. A question of Poletsky was
to know if a thin Hartogs figure has always a neighborhood
imbeddable in a Stein space. In \cite{chirka}, Chirka and
Ivashkovitch gave a counter-example arising  in an open complex
manifold. We begin by answering the question of the existence of
such figures in compact manifolds by giving a
counter-example arising in $P_2(\mathbb{C})$.\\
This two problems are related because a thin Hartogs figure being
topologically equivalent to the unit disc, it has a neighborhood
homeomorphic to the unit ball of $\mathbb{C}^2$. More, by
smoothing the thin Hartogs figure we constructed in
$P_2(\mathbb{C})$, we obtain  a smooth and closed disc with
boundary $\overline D \subset P_2(\mathbb{C})$ such that any
holomorphic function defined on any of its open neighborhoods is
constant. This disc admitting neighborhoods diffeormorphic to the
unit ball of $\mathbb{C}^2$, we obtain a negative answer to
Ivashkovitch's question.\\
For any algebraic curve $C$, $P_2(\mathbb{C}) \backslash C$ is
Stein and cannot contain $\overline D$. Thus, $\overline D$
intersects all algebraic curves of $P_2(\mathbb{C})$. Real
surfaces of $P_2(\mathbb{C})$ having this last property have been
constructed by B. Fabre \cite{fabre} then by Nemirovski
\cite{nemiro3}. Those examples are in some sense contrary to our's
because they admit Stein neighborhoods.\\
Finally, let $B$ be an open  neighborhood of $\overline D$
diffeomorphic to the unit ball and with smooth boundary $\partial
B$. Then, by combining our construction, the result of
\cite{nemiro2} and the Plemedj decomposition in $P_2(\mathbb{C})$,
we prove that $\partial B$ is an example of a smooth
CR-hypersurface of $P_2(\mathbb{C})$, diffeormorphic to the unit
sphere of $\mathbb{C}^2$ such that all continuous CR functions
defined on $\partial B$ and all holomorphic functions defined on
any connected component of $P_2(\mathbb{C})\backslash \partial B$
are constant.

\section{Envelops of holomorphy of open sets in projective space}
The study of the envelops of holomorphy has been treated by
authors as Fujita \cite{fujita},\cite{fujita2}, Takeuchi
\cite{takeuchi}, Kiselman \cite{kiselman} or Ueda \cite{ueda}. Let
us recall some well known results.

\begin{defi}
Let $X$ be a complex manifold. A domain over $X$ is a connected
complex manifold $W$ equipped with a locally biholomorphic map
$\pi: U \rightarrow X$. We say that a domain $(U,\pi)$ contains
another domain $(V,\Pi)$ if there is a map $j: V \rightarrow U$
respecting the projections, $\pi  \circ j = \Pi \circ j$. The {\it
envelope of holomorphy} $(\widetilde U, \widetilde \pi)$ of a
domain $(U,\pi)$ over $X$ is the maximal domain over $X$
containing $(U, \pi)$ such that every holomorphic function in $U$
extends holomorphically to $\widetilde U$.
\end{defi}

Let us consider the complex projective space $P_n(\mathbb{C})$ as the quotient of
 $\mathbb{C}^{n+1} \backslash \{0\}$ by the action of $\mathbb{C}^*$.
  Holomorphic functions in domains over $P_n(\mathbb{C})$ correspond to holomorphic functions
over $\mathbb{C}^{n+1}$ constant on the lines passing through the origin.
The analytic continuation in $\mathbb{C}^{n+1}$ preserves
this property because it can be represented by the differential
equation $$\sum_{j=1}^{n+1} z_j \frac{\partial f}{\partial z_j}=0.$$
Hence, envelopes over $P_n(\mathbb{C})$ correspond to envelopes over $\mathbb{C}^{n+1}$.

\begin{prop}
Let $U$ be a domain over $P_n(\mathbb{C})$, we have two cases, if the envelope over $\mathbb{C}^{n+1}$
 contains the origin, then all holomorphic functions on $U$ are constant and the envelope over
 $P_n(\mathbb{C})$ is the entire space. Otherwise,
the envelope over $\mathbb{C}^{n+1}$ is a Stein manifold.
\end{prop}
\begin{prop}
The envelope of holomorphy of a domain over $P_n(\mathbb{C})$
is either a Stein manifold or coincides with the entire $P_n(\mathbb{C})$.
 Equivalently, the envelope is Stein if and only
if there exists non-constant holomorphic functions on the domain.
\end{prop}

\section{Continuity principle}

Let $X$ be a complex manifold,  an {\it analytic disc of $X$} is a
continuous map $A: \overline \Delta \rightarrow X$ which is
holomorphic on $\Delta$. The {\it boundary} $\partial A$ of the
analytic disc $A$ is by definition the restriction of $A$ to the
unit circle $S^1=\partial \Delta$. A family of discs $\{A_t\}_{t
\in [0,1]}$ is called {\it continuous} if the map $\widetilde A:
[0,1]\times \overline \Delta \rightarrow X$ defined by $\widetilde
A(t,w)=A_t(w)$ is continuous. Let us  recall the following well
known continuity principle (see \cite{chabat}):

\begin{prop}[Behnke-Sommer]
Let $\{A_t\}_{t \in [0,1]}$ be a continuous family of analytic
discs of a complex manifold $X$. Let $U \subset X$ be an open set
and $f: U \rightarrow \mathbb{C}$ be a holomorphic function.
Suppose that $U$ verifies the following:
\begin{enumerate}
\item $A_0 \subset U$. \item For any $t \in [0,1]$, the boundary
$\partial A_t \subset U$.
\end{enumerate}
Then for any $t \in [0,1]$, their exists a neighborhood $U_t$ of
the disc $A_t$ such that $f$ extends holomorphically to $U_t$.
\end{prop}

\section{Construction of the thin Hartogs figure}}

For any point $z \in \mathbb{C}^3$, let $L_z$ be the complex line
passing through $z$ and the origin, this line defines a point
$\widetilde L_z$ in $P_2(\mathbb{C})$.  Let $\Phi:
\mathbb{C}^3\backslash \{0\} \rightarrow P_2(\mathbb{C})$ be the
map defined by $\Phi(z)=\widetilde L_z$. If $\{A_t\}_{t\in [0,1]}$
is a smooth family of closed analytic discs properly imbedded in
$\mathbb{C}^3$, such that $A_1(0)=0$, then the smooth family of
analytic discs $\{\Phi \circ A_t\}_{t\in [0,1[}$ is well defined.

\begin{prop}
Let $\cal W$ be the thin Hartogs figure $\Delta \times \{0\}\cup
S^1\times [0,1] \subset \mathbb{C}^2$. Their exists two complex
lines $L_1$ and $L_2$ of $\mathbb{C}^3$ and a continuous (even
smooth)  family $\{A_t\}_{t \in [0,1]}$ of closed analytic discs
of $\mathbb{C}^3$ such that the family of analytic discs $\{\Phi
\circ A_t\}_{t \in [0,1[}$ is
\begin{enumerate} \item continuous and properly embedded in
$P_2(\mathbb{C})$. \item For $0 \leq t_1 < t_2 < 1$ the discs
$\Phi \circ A_{t_1}$ and $\Phi \circ A_{t_2}$ intersects only at
the points $\widetilde L_1$ and $\widetilde L_2$. \item For any $t
\in [0,1[$ the disc $A_t$ is transversal to $L_1$ and to $L_2$.
\item The restriction of the map $\widetilde {\Phi \circ A}$
defined by $\widetilde {\Phi \circ A}(w,t)=\Phi\circ A_t(w)$ to
 $\cal W$ is a continuous (even smooth) proper imbedding of $\cal W$ into $P_2(\mathbb{C})$
 ($\widetilde {\Phi \circ A}({\cal W})$ is a thin hartogs figure of $P_2(\mathbb{C})$).
\end{enumerate}
\end{prop}
\begin{demo}
Let $P(z_1,z_2,z_3): \mathbb{C}^3 \rightarrow \mathbb{C}$ be a
generic polynomial of degree $2$ such that the complex
hypersurface $H=\{P(z)=0\}$ is a smooth and generic quadric which
contains the origin. Thus,  $H$ contains only two complex lines
$L_1$ and $L_2$ passing through the origin. According to the
Bezout theorem, for any point $z \in (H\backslash (L_1 \cup
L_2))$, the line $L_z$ intersects $H$ only at the point $z$ and at
the origin. Then, the restriction of the map $\Phi$ on the Zarisky
open set $H\backslash (L_1 \cup L_2)$ is open, one to one and
holomorphic (it defines a biholomorphism on its image). Let $F:
\mathbb{C}^3 \rightarrow \mathbb{C}$ be a holomorphic submersion
and note $F_c$  the smooth hypersurface $F_c=\{F(z)=c\}$. Suppose
$F$ is chosen such that $F_0$ is transversal at the origin to
$L_1$, $L_2$ and $H$. Then $F_0$ intersects $L_1$ and $L_2$ only
at the origin and intersects $H$ on a smooth curve $S_0=H \cap
F_0$. Let us note, for any $c \in \mathbb{C}$, $S_c=H \cap F_c$.
Then their exists a small neighborhood $V$ of the origin in
$\mathbb{C}$ such that $\{S_c\}_{c \in V}$ is a smooth family of
complex curves of $H$ transversal to the lines $L_1$ and $L_2$. If
$V$ is taken small enough, their exists $\epsilon >0$ such that,
for any $c \in V$ the ball $B(0,\epsilon) \subset \mathbb{C}^3$
intersects $S_c$ on an analytic disc $B_c$. One can always choose
the parametrization of the discs $B_c$ and an imbedding $\phi$ of
the set $[0,1]$ in $V$ with $\phi(1)=0$ such that the family of
disc $\{A_t\}_{t\in[0,1]}=\{B_{\phi^{-1}(t)}\}_{t \in [0,1]}$ is a
smooth and properly imbedded family of analytic discs of
$\mathbb{C}^3$ with $A_1(0)=0$. By the transversality assumption
(for $V$ chosen small enough) , this family of analytic discs
verify the Lemma.
\end{demo}

{\bf Remark}. By exploding $P_2(\mathbb{C})$ at the points
$\widetilde L_1$ and $\widetilde L_2$ (let us denote $\widetilde
P_2(\mathbb{C})$ this manifold) our construction gives an
imbedding of the family $\{A_t\}_{t \in [0,1[}$ in $\widetilde
P_2(\mathbb{C})$.







\begin{theo} Let $\{A_t\}_{t \in [0,1]}$ be the smooth family of
analytic disc constructed in the previous proposition and ${\cal
H}=\widetilde {\Phi \circ A}({\cal W})$ the corresponding thin
Hartogs figure of $P_2(\mathbb{C})$. Then  any holomorphic
function defined in a connected neighborhood of $\cal H$ is
 constant.Thus no neighborhood of ${\cal H}$ can be embedded in a
 Stein space.
 \end{theo}
\begin{demo}
Let $U$ be an open and connected neighborhood of $\cal H$ in
$P_2(\mathbb{C})$ and $f$ be a holomorphic function defined on
$U$. Let us note $\widehat U$ and $\widehat f$ the corresponding
open set and holomorphic function of $\mathbb{C}^3 \backslash
\{0\}$. Then, by construction, $\widehat U$ contains an open
neighborhood of $A({\cal W})$. According to the continuity
principle, $\widehat f$ extends holomorphically in a neighborhood
of the disc $A_1(\overline \Delta)$. So the envelop of holomorphy
of $\widehat U$ over $\mathbb{C}^3$ contains the origin and
according to proposition $1$, $f$ has to be constant.
\end{demo}

Let us note $\Delta(r)=\{ z \in \mathbb{C}; |z| \leq r \}$ with $r
\in [0,1]$, $S^1(r)$ its boundary and let us define
$$\overline D=A_1(\overline \Delta({1/2})) \cup_{t\in [0,1[} A_t(S^1({1/2+2^{\frac{1}{t-1}}})).$$
Then $\overline D$ is a smooth disc with boundary and as for the
precedent proposition, any holomorphic function defined in any of
its neighborhoods has to extend to a domain over $\mathbb{C}^3$
which contains the origin and so, has to be constant. Moreover,
for any compact complex curve  $C \subset P_2(\mathbb{C})$, the
open set $P_2(\mathbb{C}) \backslash C$ is pseudoconvex and Stein.
As there exists non-constant holomorphic functions in Stein
manifolds, the disc $\overline D$ is not included in
$P_2(\mathbb{C})\backslash C$. Thus $D$ intersects $C$. We have
obtained:
\begin{coro}
Their exists a (non analytic) closed and smooth disc with boundary
$\overline D \subset P_2(\mathbb{C})$  such that any holomorphic
function defined on its neighborhood is constant. Consequently,
$\overline D$ intersects any algebraic curves of
$P_2(\mathbb{C})$.
\end{coro}

The disc with boundary $\overline D$ being smooth, it has an open
neighborhood $B$ diffeomorphic to the unit ball of $\mathbb{C}^2$.

\begin{coro}
Their exists a domain $B \subset P_2(\mathbb{C})$, diffeomorphic
to the unit ball of $\mathbb{C}^2$ such that any holomorphic
function defined on it is constant.
\end{coro}

Moreover, $\partial B$ is a smooth hypersurface dividing
$P_2(\mathbb{C})$ into two domains $B$ and
$C_B=P_2(\mathbb{C})\backslash \overline B$. The domain $B$ being
topologically trivial, its complementary $C_B$ has the same second
holomology group than $P_2(\mathbb{C})$, in particular, $C_B$
contains a  non contractible real $2$-sphere. According to
\cite{nemiro2}, all holomorphic functions defined on $C_B$ are
constant. It is well known (see \cite{sarkis} for an example of a
proof) that the Plemedj decomposition of CR functions is available
in $P_2(\mathbb{C})$. So any continuous CR function $f$ defined on
$\partial B$, can be decomposed as $f=f^+-f^-$ with $f^+$ and
$f^-$ the boundary values (in the current sense) of holomorphic
functions defined respectively on $B$ and $C_B$. In our case, this
two holomorphic functions have to be constant, so $f$ is constant.
\begin{coro}
The boundary $\partial B$ of the previously constructed domain $B
\subset P_2(\mathbb{C})$  is a smooth hypersurface dividing
$P_2(\mathbb{C}) \backslash \partial B$  into two domains $B$ and
$C_B$ which verify the following properties:
\begin{enumerate}
\item All holomorphic functions defined on  $B$ are constant.
\item All holomorphic functions defined on $C_B$ are constant.
\item The CR hypersurface $\partial B \subset P_2(\mathbb{C})$ is
diffeomorphic to the unit sphere $S^3$ of $\mathbb{C}^2$ and all
continuous CR functions defined on $\partial B$ are constant.
\end{enumerate}
\end{coro}

 \noindent e-mail: {\tt sarkis@math.univ-lille1.fr}\par\noindent
\par
\medskip

\begin{thebibliography}{9}

\bibitem{chabat}
Chabat B. {\em Introduction \`a l'analyse complexe}, Tome 2, Fonctions de plusieurs variables, Editions Mir, p. 215.

\bibitem{chirka}
Chirka E. and Ivashkovich S., {\em On nonimbeddability of Hartogs
figures into complex manifolds}, Arxiv.math.CV/0404290, 16 April
2004.
\bibitem{fabre}
Fabre B. {\em Sur l'intersection d'une surface de Riemann avec des
hypersurfaces algébriques}, C.R.Acad. Sci. Paris Sér I Math., {\bf
322} (1996), no. 4, 371-376.

\bibitem{fujita}
Fujita R., {\em Domaines sans point critiques int\'erieur sur l'espace
projectif complexe}, J. Math. Soc. Japan, {\bf 15} (1963), no. 4,
443-473.
\bibitem{fujita2}
Fujita R., {\em Domaines sans point critiques int\'erieur sur l'espace
  produit}, J. Math. Kyoto Univ., {\bf 4} (1965), no. 3, 493-514.


\bibitem{kiselman}
Kiselman C.O., {\em On entire functions of exponential type and
  indicators of analytic functions}, Acta Math., {\bf 117} (1967), 1-35.

\bibitem{nemiro3}
Nemirovski S. Yu.,{\em Stein domains on algebraic varieties}, Mat.
Zametki {\bf 60} (1996), 295-298; English transl., Math. Notes
{\bf 60} (1996), 218-221.

\bibitem{nemiro2}
Nemirovski S. Yu.,{\em Holomorphic functions and embedded real
surfaces"}, Mat. Zametki {\bf 63} (1998), 599-606; English
transl., Math. Notes {\bf 63} (1998), 527-532.

\bibitem{nemiro}
Nemirovski S. Yu.,{\em Complex analysis and differential topology
on complex surfaces}, Russian Math. Surveys {\bf 54} (1999),
729-752.

\bibitem{sarkis}
Sarkis F., {\em CR-meromorphic extension and the nonembeddability
of the Andreotti-Rossi CR structure in the projective space},
Internat. J. Math., {\bf 10} (1999), no. 7, 897-915.

\bibitem{sarkis2}
Sarkis F., {\em Hartogs-Bochner type theorem in projective space},
Ark. Math., {\bf 41} (2003), no. 1, 151-163.

\bibitem{takeuchi}
Takeuchi A., {\em Domaines pseudoconvexes infinis et la riemannienne
  dans un espace projectif}, J. Math. Soc. Japan, {\bf 16} (1964), 159-181.


\bibitem{ueda}
Ueda T., {\em Pseudoconvex domains over Grassmann manifolds}, J.
Math of Kyoto Univ, {\bf 20} (1980), 391-394.


\end{thebibliography}
\end{document}